%% file: main.tex
\documentclass[journal, twoside]{IEEEtran}

\usepackage{import}
\subimport{praemble/}{hft-paper-praemble-main.tex}

\usepackage[style=ieee, backend=biber, isbn=false,maxbibnames=99]{biblatex}
  
\graphicspath{{figures/}} 
\usepackage{atbegshi}

\begin{document}

    \title{An Explicit Higher-Order Dual Basis for a Multiplicatively Calderón Preconditioned \\ Electric Field Integral Equation}

	%

    \author{Bernd Hofmann,~\IEEEmembership{Member,~IEEE,}
	        Thomas F. Eibert,~\IEEEmembership{Senior Member,~IEEE,\\}
	        Francesco P. Andriulli,~\IEEEmembership{Fellow,~IEEE,}
			and~Simon B. Adrian,~\IEEEmembership{Senior Member,~IEEE}
	\thanks{This work was funded by the Deutsche Forschungsgemeinschaft (DFG, German Research Foundation) – 504345461.}%
	\thanks{B. Hofmann is with the Department of Electrical Engineering, Stanford University, Stanford, CA 94305, USA (e-mail: bernd.hofmann@stanford.edu).}
    \thanks{T. F. Eibert is with the Department of Electrical Engineering, School of Computation, Information and Technology, Technical University of Munich, 80290 Munich, Germany.}
	\thanks{F. P. Andriulli is with the Department of Electronics and Telecommunications, Politecnico di Torino, 10129 Turin, Italy.}%
	\thanks{S. B. Adrian is with the Fakultät für Informatik und Elektrotechnik, Universität Rostock, 18059 Rostock, Germany (e-mail: simon.adrian@uni-rostock.de).}
	}

	%
	%

	\markboth{}
	{Hofmann \MakeLowercase{\textit{et al.}}: An Explicit Higher-Order Dual Basis for a Multiplicatively Calderón Preconditioned EFIE}
	%



	\maketitle

    \begin{abstract}
One of the most effective means to precondition the \ac{EFIE} discretized with \ac{RWG} functions is the multiplicative Calderón preconditioner employing \ac{BC} functions as a basis dual to the \ac{RWG} basis.
It results in a formulation that is free from the dense-discretization and the low-frequency breakdown.
To generalize the multiplicative Calderón preconditioner from the low-order \ac{BC} and \ac{RWG} basis to higher orders, we utilize B-spline-based basis functions and establish the first explicit high-order dual basis.
It can be regarded as a generalization of the \ac{BC} functions to arbitrary polynomial degrees and constitutes a fundamental building block for other approaches that rely on a dual basis.
Numerical results for the obtained preconditioner demonstrate a low and constant number of \ac{GMRES} iterations independent of the number of unknonws and the polynomial degree for canonical and realistic \ac{PEC} scatterers; a key to enable the full potential of higher-order bases.
    \end{abstract}
    \acresetall

\vspace{-0.5mm}
    \begin{IEEEkeywords}
B-splines, Calderón, EFIE, integral equations, isogeometric, low frequency, multiply connected, NURBS.
    \end{IEEEkeywords}
\vspace{-3mm}

    %
    \IEEEpeerreviewmaketitle

    \section{Introduction}

\IEEEPARstart{G}{alerkin} discretizations of the \ac{EFIE}~\cite{maueZurFormulierungAllgemeinen1949} have been shown to be well suited to numerically solve time-harmonic electromagnetic radiation and scattering problems involving \ac{PEC} objects.
The Galerkin scheme guarantees convergence to the exact solution of the \ac{EFIE} when increasing the number of basis functions approximating the unknown surface current density~\cite{buffagalerkin2003}.
However, such an increase in the number of basis functions leads to the dense-discretization breakdown, that is, the condition number of the matrix in the \ac{LSE} to be solved increases.
An effective means to overcome this breakdown is Calderón preconditioning\footnote{Note that Calderón preconditioning also partly cures the low-frequency breakdown. However, to fully overcome the latter a combination with quasi-Helmholtz projectors~\cite{andriulliWellConditionedElectricField2013,hofmannLowFrequencyStabilizationBSpline2024,hofmannExcitationAwareSelfAdaptiveFrequency2023,hofmannLowFrequencyStableExcitation2023} is needed (and possible).}: leveraging on the self-regularizing property of the \ac{EFIE}, the Calderón identity, the \ac{EFIE} operator can be applied to itself to yield a second-kind integral operator composed of an identity operator and a compact pertubation~\cite{steinbachConstructionEfficientPreconditioners1998,hiptmairOperatorPreconditioning2006,buffaDualFiniteElement2007,andriulliMultiplicativeCalderonPreconditioner2008,antoineIntroductionOperatorPreconditioning2021}.
A purely multiplicative Calderón preconditioner of this kind has been established in~\cite{buffaDualFiniteElement2007,andriulliMultiplicativeCalderonPreconditioner2008,changQuadrilateralBarycentricBasis2013} for the \ac{RWG} (primal) basis~\cite{raoElectromagneticScatteringSurfaces1982} by employing \ac{BC} functions as a dual basis.
Such a dual basis is needed to link the two \ac{EFIE} operators in a way that preserves the properties of the analytical Calderón identity in its discrete counterpart.

A drawback of this approach is that the \ac{RWG} and \ac{BC} functions are lowest order, which limits the convergence in accuracy.
For the \ac{RWG} basis, different generalizations to higher polynomial degrees have been proposed (e.g.,~\cite{petersonPARAMETRICMAPPINGVECTOR1995,gragliaHigherOrderInterpolatory1997,petersonMappedVectorBasis2006}), resulting in schemes that yield the same accuracy with fewer unknowns and thus, reduce computation and memory requirements. 
However, for the \ac{BC} functions, no direct generalization to higher orders is available:
The approach of~\cite{valdesHighorderDivQuasi2011} constructs a dual basis numerically based on \acp{SVD} and local orthogonalizations on a barycentrically refined mesh;
the numerical results indicate that the condition numbers obtained using the local orthogonalization are higher than the ones obtained with the---in practice computationally prohibitive---global orthogonalization, thus indicating that the approach is suboptimal. 
The approach of~\cite{bourhisNovelCalderonPreconditioning2024,bourhisPreconditioningStrategiesConformal2025} uses projections from spaces with higher polynomial degree than the degree of the primal basis.
The self-dual basis in~\cite{liSubdivisionBasedIsogeometric2016} is limited to smooth geometries and fixed to cubic polynomials.

Instead, we are interested in constructing an explicitly defined higher-order dual basis.
This enables a direct transfer~of the tools known from the lowest order to higher orders.
To this end, we consider the divergence-conforming B-spline-based basis functions proposed in~\cite{hofmannLowFrequencyStabilizationBSpline2024} and similar to~\cite{veigaMathematicalAnalysisVariational2014,simpsonIsogeometricBoundaryElement2018,dolzIsogeometricBoundaryElements2019} as the primal basis on curvilinear quadrilateral patches. 
In the lowest order, they correspond to \ac{RWG} functions on quadrilateral subdomains of the patches (also called rooftop~\cite{petersonMappedVectorBasis2006} or Raviart-Thomas functions~\cite{raviartMixedFiniteElement1977}) corresponding to a virtual Cartesian splitting of the patches.
In contrast to Lagrange-based functions, they have local support and allow for an arbitrary number of basis functions per patch, independent of the polynomial degree.
Originally, the B-splines were introduced as basis functions in the context of structural mechanics~\cite{hughesIsogeometricAnalysisCAD2005}, and the term isogeometric analysis was coined as the underlying geometry can also be defined using B-splines, commonly via \ac{NURBS}.
The latter often allow for exact representation of geometries (which are often in the first place defined via \ac{NURBS}) in contrast to geometrical discretization errors introduced by classical meshes. 

In this work, we derive the first explicit higher-order dual basis for a primal divergence-conforming B-spline basis and construct a higher-order Calderón preconditioner.
To this end, we introduce a refined parametric space---akin to the barycentric refinement used in the case of the \ac{RWG} space.
The latter is defined precisely so that the dual basis can be obtained as a local linear combination of refined basis functions, with scalar coefficients that are trivial to compute.
In the lowest-order case, the obtained dual basis is identical to the \ac{BC} functions and for the higher orders it can, hence, be regarded as a generalization of the latter to arbitrary polynomial degrees.
By introducing a non-uniform knot vector for the dual basis, we obtain a mixed Gram matrix with a condition number independent of the number of unknonws.
The presented approach is applicable to open and closed, simply- and multiply-connected geometries, which may be smooth or include edges and corners.
Numerical results for both canonical and realistic geometries corroborate the effectiveness of this strategy.
Note that preliminary results have been presented in~\cite{hofmannICEAA2024,hofmannKleinheubach2024}.

This article is organized as follows: 
Section~II introduces background material about the B-spline-based discretization of the \ac{EFIE} and fixes the notation. 
In Section~III, the higher-order primal basis is defined, the corresponding higher-order dual basis is derived for arbitrary polynomial degrees on closed and open geometries, and all relevant properties are discussed.
The efficiency of the resulting Calderón preconditioner is demonstrated for the scattering from canonical and realistic objects in Section~IV, followed by a conclusion.

    \section{B-Spline Based Discretization of the EFIE}

Let $\Gamma \subset \mathbb{R}^3$ denote the surface of a \ac{PEC} scatterer immersed in a homogeneous background medium characterized by permittivity $\varepsilon$ and permeability $\mu$.
The surface $\Gamma$ is assumed to be a bounded, two-dimensional Lipschitz manifold, which can be simply or multiply connected, open or closed, equipped with a unit normal vector $\veg n$.
When excited by a time-harmonic incident field $(\veg e^\mr{ex}, \veg h^\mr{ex})$, a surface current density $\veg j$ is induced, which is governed by the \ac{EFIE}~\cite{mautzHFieldEFieldCombinedField1978} 
\begin{equation}
    \vecop T_{\!\!k} \veg j = \veg n \times \veg e^\mr{ex} \,,
   \label{efie}
\end{equation} 
where the operator $\vecop T_{\!\!k} \veg j = \jm k \TAop \veg j + \jm k^{-1} \TPop \veg j$ is composed of the vector potential operator
\begin{equation}
    \TAop \veg j =  \veg n \times \iint _\Gamma G_k(\veg r, \veg r') \, \veg j(\veg r') \,\dd S(\veg r')
\end{equation}
and the scalar potential operator
\begin{equation}
    \TPop \veg j = \veg n \times \nabla \! \iint _\Gamma G_k(\veg r, \veg r') \,\nabla_\Gamma \cdot \veg j(\veg r') \,\dd S(\veg r') \,.
    \label{TP}
\end{equation}
Here, $G_k(\veg r, \veg r') = \e^{-\jm k |\veg r - \veg r'|} \big/ ({4\uppi |\veg r - \veg r'|})$ is the free-space Green's function of the scalar Helmholtz equation, $k=\omega \sqrt{\mu\varepsilon}$ denotes the wavenumber, $\omega$ the angular frequency, and \mbox{$\jm^2=-1$} the imaginary unit.
Moreover, the product with the wave impedance is absorbed into the current density $\veg j$, and a suppressed time convention of $\e^{\,\jm \omega t}$ is assumed.
Once $\veg j$ has been determined, the corresponding scattered or radiated fields follow from it.

    \subsection{Galerkin Discretization}

To solve the \ac{EFIE} in~\eqref{efie} for $\veg j$, the latter is approximated by an expansion in basis functions $\veg f_{\!\!n}$, i.e., $\veg j \approx \sum_{n=1}^N {[\vec j]}_n \veg f_{\!\!n}(\veg r)$, where $\vec j \in \mathbb{C}^N$ contains the unknown expansion coefficients (note that we employ sans-serif letters for discretized quantities). 
The basis functions are defined on the surface $\Gamma$, which we describe by a union of $N_A$ patches $\Gamma_a$, whose pairwise intersections are either empty, a vertex, or an edge, that is, $\Gamma = \bigcup_{a=1}^{N_A} \Gamma_a$, where each $\Gamma_a$ is described by a mapping $\veg{s}_a(u,v): \hat{\Gamma}_a \mapsto \Gamma_a$ from the parametric space $\hat{\Gamma}_a = (u, v) \in {[0,1]}^2$ to the physical space $\Gamma_a \subset \mathbb{R}^3$.
Consequently, we employ basis functions $\fhatn(u,v)$ defined on the unit square $(u,v) \in [0,1]^2$ and mapped to rectangular curvilinear patches as (see, e.g.,~\cite{petersonMappedVectorBasis2006})
\begin{equation}
    \veg f_{\!\!n}(\veg r) = \cfrac{\veg J(u,v)}{D(u,v)} \,\fhatn(u,v) 
    \label{fnMap}
\end{equation}
with $\veg r = \veg{s}_a(u,v)$, the Jacobian matrix $\veg J \in \mathbb{R}^{3 \times 2}$ of the mapping $\veg{s}_a$, and the corresponding generalized Jacobian determinant $D$. 
Note that throughout this article, we denote with a hat the parametric counterparts $\hat{\veg g}$ or $\hat g$ to a vector $\veg g$ or a scalar $g$ in the physical domain.

Inserting the expansion of $\veg j$ into~\eqref{efie} and testing with the same $\veg f_{\!\!n}$, the \ac{LSE}
\begin{equation}
    \vec{T}_{\!k}\vec{j} = (\jm k \TA + \jm k^{-1} \TP) \vec j  = \vec{e}^\mr{ex} \,,
    \label{efieLSE}
\end{equation}
is obtained, where the matrix $\TA \in \mathbb{C}^{N\times N}$ consists of the entries 
\begin{equation}
    \left[\TA\right]_{mn} = \iint_\Gamma \veg f_{\!\! m} \cdot  \iint_{\Gamma} G_k(\veg r, \veg r')  \veg f_{\!\! n} \, \dd S(\veg r') \,\dd S(\veg r)\,              
    \label{singular}
\end{equation}
the matrix $\TP \in \mathbb{C}^{N\times N}$ consists of the entries 
\begin{equation}
    \left[\TP\right]_{mn} = - \!\! \iint_{\Gamma}\! \nabla_{\!\Gamma} \cdot \veg f_{\!\! m} \! \iint_{\Gamma} G_k(\veg r, \veg r')  \nabla_{\!\Gamma} \cdot \veg f_{\!\! n} \, \dd S(\veg r') \dd S(\veg r)\,,
    \label{hsingular}
\end{equation}
and the \ac{RHS} vector $\vec e^\mr{ex} \in \mathbb{C}^N$ consists of the entries ${\left[\vec e^\mr{ex}\right]}_m = \iint_\Gamma \veg f_{\!\!m} \cdot \veg e^\mr{ex}\,\dd S(\veg r) $.

    \subsection{Multiplicative Calderón Preconditioning}

Leveraging on the Calderón identity $\vecop T_{\!\!k}^2 = -\vecop I/4 + \vecop K^2$, where $\vecop I$ is the identity operator and $\vecop K$ a compact operator, the multiplicatively Calderón preconditioned \ac{EFIE}~\cite{buffaDualFiniteElement2007,andriulliMultiplicativeCalderonPreconditioner2008,hiptmairOperatorPreconditioning2006}
\begin{equation}
    \vec G^{-\T} \widetilde{\vec T}_{\!\tilde k} \vec G^{-1} \vec T_{\!k} =  \vec G^{-\T} \widetilde{\vec T}_{\!\tilde k} \vec G^{-1} \vec e^\mr{ex}\,
    \label{GTGTcalderon}
\end{equation}
can be constructed, where 
\begin{equation}
    \left[\vec G\right]_{mn}= \iint_\Gamma \veg n \times \veg f_{\!\!m} \cdot \ftil_{\!\!n}\,\dd S(\veg r)
    \label{gramDef}
\end{equation}
is a Gram matrix and 
\begin{equation}
    {[\widetilde{\vec T_{\!\tilde k}}]}_{mn} = \iint_\Gamma  \veg n \times \ftil_{\!m} \cdot \vecop T_{\!\!\tilde k} \ftil_{\!n}\,\dd S(\veg r)
\end{equation}
is a discretization of the \ac{EFIE} analogous to~\eqref{efieLSE}-\eqref{hsingular} but with a set of suitable, \emph{and yet to be defined}, dual functions $\ftil_{\!n}$.
Choices for the wavenumber $\tilde k$ are depending on the application and may be $\tilde k = k$, or complexified versions to avoid introducing interior resonances, such as $\tilde k = -\jm k$~\cite{adrianElectromagneticIntegralEquations2021,merliniMagneticCombinedField2020} or $\tilde k = k - \jm 0.4 k^{1/3}H^{2/3}$, where the latter also addresses the high-frequency breakdown with $H$ denoting the (maximum) mean curvature of $\Gamma$ \cite{adrianElectromagneticIntegralEquations2021,DarbasDiss,boubendirWellconditionedBoundaryIntegral2014}. 
The efficiency of this preconditioning strategy is rooted in the bounded spectrum of $\vecop T_{\!\!k}^2$.
Note that the inverse of the sparse matrix $\vec G$ and its transpose can be obtained efficiently and to a high accuracy, e.g., iteratively by employing an incomplete LU preconditioner.

	\section{Proposed Primal and Dual Basis}

For a combination of primal basis functions $\veg f_{\!\!n}$ and dual basis functions $\ftil_{\!n}$ to be \emph{suitable}, three properties have to be necessarily satisfied~\cite{buffaDualFiniteElement2007,andriulliMultiplicativeCalderonPreconditioner2008,coolsTimeDomainCalderon2009}: 
\renewcommand{\theenumi}{\roman{enumi}}
\begin{enumerate}
    \item both are divergence conforming,
    \item they lead to a well-conditioned Gram matrix $\vec G$ defined in~\eqref{gramDef}, and
    \item they ensure that
\end{enumerate}
\begin{equation}
    \TPd \vec G^{-1} \TP \approx \matO
    \label{vanCond}
\end{equation}
holds to machine precision with $[\widetilde{\vec T}_{\!\upPhi,\tilde k}]_{mn} = \iint_\Gamma \veg n \times \ftil_{\!\!m} \cdot \vecop T_{\!\!\upPhi,\tilde k} \ftil_{\!\!n}\,\dd S(\veg r)$ denoting the dually discretized version of the operator~\eqref{TP} analogous to~\eqref{hsingular}.
To satisfy all of the above properties, we propose the following pair of primal and dual bases.

    \subsection{Primal Basis}

\begin{figure}
	\centering
	\includegraphics[scale=1]{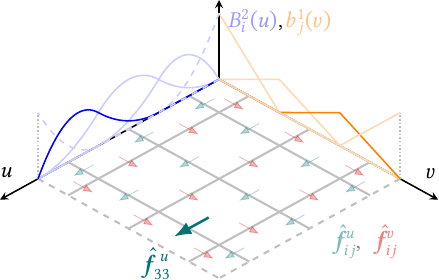}
	\caption{Primal basis functions $\fhatijd$ in the parametric domain $(u,v) \in [0,1]^2$ associated with the primal Greville mesh $\mathscr{M}_\mr{G}$ for $N_u=N_v=N_{uv}=5$ and $p_u=p_v=p_{uv}=2$.}
	\label{primal}
\end{figure}

For the primal basis of arbitrary polynomial degree, we choose the same divergence conforming basis functions as introduced in~\cite{hofmannLowFrequencyStabilizationBSpline2024}, since we can leverage the loop-star decomposition of~\cite{hofmannLowFrequencyStabilizationBSpline2024} and its specific graph properties. 
That is, for the $\fhatn$ in \eqref{fnMap}, we set $\fhatn = \fhatijd$ on a single patch with
\begin{equation}
    \fhatij^u(u,v) = B_i^{p_u}(u) b_j^{p_v-1}(v) \hat{\veg e}_u \quad \begin{array}{l} \, i = 2, 3, \dots, N_u-1 \\  j=2,3,\dots,N_v \end{array} 
    \label{fu}
\end{equation}
pointing in the $u$ unit direction $\hat{\veg e}_u$ as well as
\begin{equation}
    \fhatij^v(u,v) = b_i^{p_u-1}(u) B_j^{p_v}(v) \hat{\veg e}_v \quad \begin{array}{l} \, i = 2, 3, \dots, N_u \\  j=2, 3, \dots, N_v-1 \end{array} 
    \label{fv}
\end{equation}
pointing in the $v$ unit direction $\hat{\veg e}_v$ with the polynomial degrees $p_u\geq 1$ and $p_v \geq 1$ and an example depicted in Fig.~\ref{primal}.
The freely selectable index mappings $i = i(n)$, $j = j(n)$, and $d=d(n)\,\in \{u,v\}$ relate Cartesian and linear indices.
The $N_u$ B-splines $B_i^{p_u}(u)$ are defined on an open knot vector
\begin{equation}
    U = \{u_1, u_2, \dots, u_{M_u}\},~~ u_i\in [0,1]
\end{equation}
of length $M_u$ with $u_i \leq u_{i+1}$, and a maximum multiplicity $p_u-1$ for the internal knots.
Analogously, $N_v$ B-splines $B_j^{p_v}(v)$ are defined on the open knot vector $V = \{v_1, v_2, \dots, v_{M_v}\},v_i \in [0,1]$.
Moreover, $b_i^{p_u}$ and $b_j^{p_v}$ denote Curry-Schoenberg splines~\cite[p.~88]{carldeboorPracticalGuideSplines2001} (i.e., normalized B-splines) defined on the same knot vectors $U$ and $V$.
The choice of employing Curry-Schoenberg splines is crucial, as it facilitates the exchange of the solenoidal and non-solenoidal subspaces between primal and dual basis.
Note that via the indexing limits in~\eqref{fu} and~\eqref{fv}, we have removed basis functions which would correspond to current flowing off the boundaries of $\Gamma$ as well as basis functions that are zero due to the definition of B-splines of polynomial degree $p-1$ on an open knot vector for degree $p$.
Related to this choice, a crucial relation we will make use of in the following is that the effective number $n_u$ of splines $b_i^{p_u-1}$ (different from zero) on the knot vector $U$ is 
\begin{equation}
    n_u = N_u - 1
    \label{nu}
\end{equation}
and analogously $n_v = N_v - 1$.
At the interfaces between patches, neighboring B-splines normal to the interface need to be combined as detailed in~\cite{hofmannLowFrequencyStabilizationBSpline2024} employing the same notation.  
\begin{figure}
	\centering
	\includegraphics[scale=1]{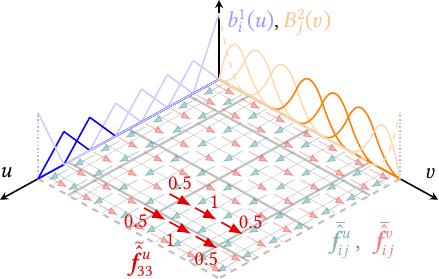}
	\caption{Refined basis functions $\ghatijd$ and example of a basis function dual to the primal basis function shown in Fig.~\ref{primal} as a linear combination of $\ghatijd$ on the refined Greville mesh $\bar{\mathscr{M}}_\mr{G}$ with $\bar N_{uv}=9$ and $p_{uv}=2$.}
	\label{dual}
\end{figure}

    \subsection{Dual Basis -- Construction}

To derive a dual basis, we utilize the fact that the basis functions in~\eqref{fu} and~\eqref{fv} can be associated with the edges of a so-called Greville mesh $\mathscr{M}_\mr{G}$~\cite{buffaIsogeometricMethodsComputational2014} as shown in Fig.~\ref{primal}.
Note that this mesh is purely virtual and serves only to visualize the underlying structure of the B-spline basis.
The mesh $\mathscr{M}_\mr{G}$ follows from the location of the Greville sites 
\begin{equation}
    \gamma_i^u =  \dfrac{1}{p_u} \sum_{j=1}^{p_u} u_{i+j}
\end{equation}
for $p_u>0$ with respect to the knot vector $U$ and analogously $\gamma_i^v$ with respect to $V$.
Crucially, also the loop-star decomposition of the primal basis 
\begin{equation}
    \vec j = \matL \vec j_{\matL} + \vec H \vec j_{\vec H}  + \matS  \vec j_\matS 
\end{equation}
into solenoidal (loop) components $\vec j_{\matL}$ and $\vec j_{\vec H}$ as well as non-solenoidal (star) components $\vec j_\matS$ via the local loop mapping matrix $\matL$, the global loop mapping matrix $\vec H$, and the star mapping matrix $\matS$ (see \cite{hofmannLowFrequencyStabilizationBSpline2024} for definition and properties), is associated with $\mathscr{M}_\mr{G}$, that is, stars are associated with the faces and loops with the vertices of $\mathscr{M}_\mr{G}$~\cite{hofmannLowFrequencyStabilizationBSpline2024}.
This motivates finding a dual basis that can be associated with $\mathscr{M}_\mr{G}$ such that the dual stars can be associated with the vertices and dual loops with the faces of $\mathscr{M}_\mr{G}$ (i.e., exchanging the roles of loops and stars), which is a key ingredient for property iii) as will be shown in Section~\ref{secProperties}.
To achieve precisely this, we propose to first introduce a refined basis exemplified in Fig.~\ref{dual}
\begin{equation}
    \ghatiju(u,v) = \bar B_i^{p_u}(u) \bar b_j^{p_v-1}(v) \hat{\veg e}_u \quad \begin{array}{l} \, i = 1, 2, \dots, \bar N_u \\  j=2,3,\dots, \bar N_v \end{array} 
\end{equation}
and
\begin{equation}
    \ghatijv(u,v) = \bar b_i^{p_u-1}(u) \bar B_j^{p_v}(v) \hat{\veg e}_v \quad \begin{array}{l} \, i = 2, 3, \dots, \bar N_u \\  j=1, 2, \dots, \bar N_v \end{array} 
\end{equation}
(which is mapped to the physical domain following~\eqref{fnMap}) where the B-splines $\bar B_i^{p_u}$ and $\bar B_j^{p_v}$ as well as the Curry-Schoenberg splines $\bar b_j^{p_v-1}$ and $\bar b_i^{p_u-1}$ are defined on the knot vectors $\bar U$ and $\bar V$.
These knot vectors correspond to an increase in the number of B-spline basis functions exactly to
\begin{equation}
    \bar N_u = 2 N_u - 1
    \label{Nref}
\end{equation}
and analogously $\bar N_v = 2 N_v - 1$ (the precise location of the knots in $\bar U$ and $\bar V$ will be addressed in Section~\ref{secGramCond}).
In doing so, we obtain a refined (purely virtual) Greville mesh $\bar{\mathscr{M}}_\mr{G}$ which can be interpreted as a generalization of the barycentrical refinement for the definition of the lowest-order \ac{BC} functions in~\cite{buffaDualFiniteElement2007,andriulliMultiplicativeCalderonPreconditioner2008}.

The refinement according to~\eqref{Nref} corresponds to the following consideration: given the primal basis functions $\fhatqqu = B_i^{p_u}(u) b_j^{p_v-1}(v) \hat{\veg e}_u$, we introduce for each $B_i^{p_u}(u)$ two refined $\bar b_i^{p_u-1}(u)$ except for the first and last $B_i^{p_u}(u)$ where only one refined $\bar b_i^{p_u-1}(u)$ is needed.
Hence, and leveraging~\eqref{nu}, we need 
\begin{equation}
    \bar n_u = \bar N_u - 1 = 2 N_u - 2
\end{equation}
refined splines.
Solving for $\bar N_u$ yields~\eqref{Nref}.
For each $b_j^{p_v-1}(v)$ of the primal basis $\fhatqqu$, we introduce three refined $\bar B_j^{p_v}(v)$ splines, where two consecutive refined splines are shared among consecutive $\fhatqqu$.
Hence, we need 
\begin{equation}
    \bar N_v = 2 n_v + 1 
\end{equation}
refined splines. 
Inserting $n_v = N_v - 1$ from~\eqref{nu}, yields again~\eqref{Nref}.
In consequence, refining as proposed in~\eqref{Nref} results in the desired amount of splines to form the dual basis.
The same logic applies to the primal basis functions $\fhatqqv$. 

\begin{figure}
	\centering
	\includegraphics[scale=1]{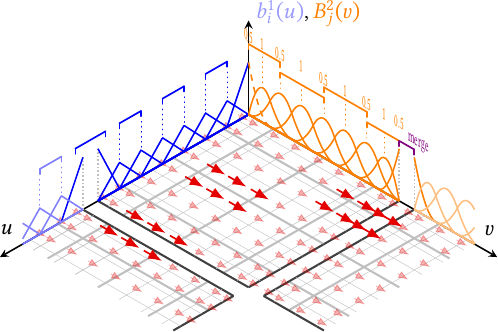}
	\caption{Forming dual basis functions from the refined basis: six refined basis functions form a dual basis function with weighting factors of 1 and 1/2. Across patch boundaries, normal continuity is enforced, and dual basis functions along the patch interface involve refined basis functions from both patches.}
	\label{dualAll}
\end{figure}

The actual dual basis functions $\ftilij$ are then formed by linearly combining the refined basis functions $\ghatijd$ (in the parametric domain) as\footnote{Note that the dual functions for \smash{$\fhatij^u$} point in the $v$-direction and the dual functions \smash{$\fhatqqv$} point in the $u$-direction.} 
\begin{equation}
    \fhattildeuij = \sum_{s,t} \alpha_{st}^u \, \ghatstv(u,v) 
\end{equation}
and
\begin{equation}
    \fhattildevij = \sum_{s,t} \alpha_{st}^v \, \ghatstu(u,v) 
\end{equation}
as exemplified in Fig.~\ref{dual}, where all coefficients $\alpha_{ij}^d$ are explicitly available.
More precisely, each dual function is fundamentally formed by 6 adjacent refined basis functions with $\alpha_{st}^{u/v} \in \{1, 1/2\}$ (with special cases for interfaces between patches addressed later on). 
The refined basis functions involved are two consecutive refined functions transverse to the basis function direction (none shared among the basis functions), and 3 consecutive refined basis functions along the basis function direction with one refined basis function being shared (which corresponds to and explains the refinement strategy~\eqref{Nref}).
This is concretized in the example shown in Fig.~\ref{dualAll} for $p_u=p_v=2$, $N_u=N_v=5$, and $\bar N_u= \bar N_v=9$.
Due to the coefficients $\alpha_{st}^{u/v} \in \{1, 1/2\}$, four dual basis functions running along the edges of a primal Greville mesh face can form a solenoidal (divergence-free) function: it can be verified that the refined basis functions actually form nine loops around the vertices of $\bar{\mathscr{M}}_\mr{G}$ with weights of 2, 1, and $1/2$.
At the patch boundaries, two geometry scenarios have to be handled: closed and open surfaces.

    \subsubsection{Closed Surfaces} \label{secClosedSurfaces}

\begin{figure}
	\centering
	\includegraphics[scale=1]{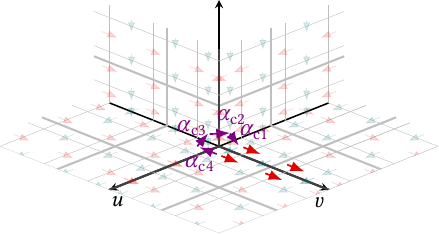}
	\caption{Forming dual basis functions at a patch corner where $N_\mr{c}$ patches meet. The weighting coefficients $\alpha_{\mr{c}i}$ ensure that solenoidal dual functions can be constructed around the faces of the primal Greville mesh.}
	\label{dualCorner}
\end{figure}

For closed surfaces, each patch shares an edge with four other patches. 
Moreover, at each patch corner, at least two patches meet.
To form suitable dual basis functions for this setup, we distinguish (next to dual basis functions completely on the interior) three different cases that need to be handled, with two of them exemplified in Fig.~\ref{dualAll}: 
a)~dual basis functions normal to a patch interface have to be complemented by the corresponding refined basis functions~$\ghatijd$ on the adjacent patch in order to ensure continuity of the normal component across the interface;
b)~dual basis functions running along an interface have to be formed by combining the refined basis functions on both sides of the interface;
c)~as illustrated in Fig.~\ref{dualCorner}, at the corner points where multiple patches meet, all refined basis functions normal to the edges incident to the corner contribute to the dual functions surrounding the corner with appropriately chosen weighting coefficients.

For case~c), the weighting coefficients must be chosen such that the dual basis functions formed around any face of the primal Greville mesh $\mathscr{M}_\mr{G}$ can be combined to yield solenoidal functions, which is essential for property iii) to hold.
Hence, and similar to the lowest order case in~\cite{buffaDualFiniteElement2007,changQuadrilateralBarycentricBasis2013}, for an interior corner where $N_\mr{c}$ patches meet, the weighting coefficients have to be chosen as 
\begin{equation}
    \alpha_{\mr{c}i} = 1 - \dfrac{2i}{N_\mr{c}} \qquad \text{for}\qquad i = 1, \dots, N_\mr{c} - 1
\end{equation}
where the index $i$ enumerates the refined basis functions around the corner.

Note that the above construction is valid independent of the polynomial degree. 
That this yields a valid dual basis is owed to the Curry-Schoenberg normalization of the involved splines, the partition of unity property of the B-splines, and the underlying structure which allows for the association of the primal as well as the dual basis functions to the edges of the Greville mesh $\mathscr{M}_\mr{G}$.

    \subsubsection{Open Surfaces}

\begin{figure}
	\centering
	\includegraphics[scale=1]{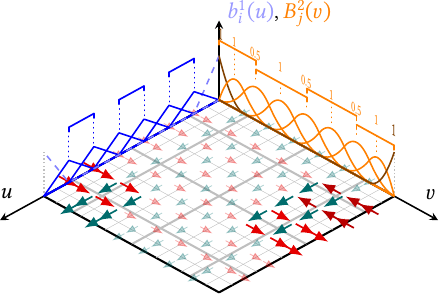}
	\caption{Forming dual basis functions on a patch touching the boundaries of $\Gamma$ by including half basis functions adjacent to the boundary with a weighting factor of $1$ instead of $1/2$. Examples of solenoidal functions formed by three dual functions touching an edge and by two dual functions touching a corner.}
	\label{dualOpen}
\end{figure}

At the boundaries of the surface $\Gamma$, present only for open surfaces, the dual basis functions can still be defined following the construction principles described above.
However, they must be modified to account for the absence of neighboring patches across the boundary.
Specifically, the dual basis functions adjacent to the boundary incorporate half basis functions of $\ghatijd$ with a weight of $1$ instead of $1/2$, as illustrated in Fig.~\ref{dualOpen}.
This modification preserves the solenoidal nature when combining the dual basis functions around the faces of the Greville mesh (which is essential for property iii) to hold also for open geometries): due to the Curry-Schoenberg normalization, the weights can be interpreted as the amount of charge flowing in and out of each Greville mesh face.
Hence, it can be verified that the net charge is indeed zero when combining dual functions at edges and corners of $\Gamma$ to a solenoidal one.
The only two possible cases for forming solenoidal functions at boundaries of $\Gamma$ are illustrated in Fig.~\ref{dualOpen}, for which the solenoidal nature is found to hold.
For all other interfaces between patches not located at the boundary of $\Gamma$, the same construction scheme as defined for closed surfaces in the previous Section~\ref{secClosedSurfaces} applies.

    \subsection{Dual Basis -- Refinement Strategy} \label{secGramCond}

Having fundamentally defined the dual basis, the precise refinement strategy (i.e., the location of the knots in $\bar{U} $ and $\bar{V}$) for the knot vectors needs to be addressed carefully: A naive choice for the refined knot vectors $\bar{U} $ and $\bar{V}$ would be a uniform refinement of the primal knot vectors $U$ and $V$, that is, the refined knots are chosen equidistant.
However, since the Greville sites of the primal basis are in general not uniformly distributed for $p_u, p_v >1$---in particular near the boundaries of the parametric domain, where they cluster due to the open knot vector construction---a uniformly refined dual basis would exhibit a growing mismatch with the primal basis as $N_u$ and $N_v$ increase.
This mismatch causes the overlap between the supports of dual and primal basis functions to become increasingly uneven and eventually leads to a Gram matrix $\vec G$ whose condition number grows with the number of unknowns.
To prevent such behavior, we propose a non-uniform refinement of the knot vectors that seeks to align the Greville sites of the dual basis with those of the primal basis.

\begin{figure}
	\centering
	\includegraphics[scale=1]{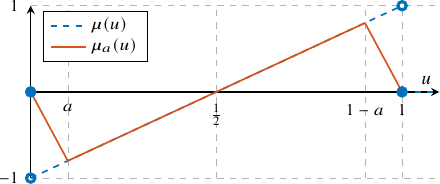}
	\caption{Functions $\mu(u)$ and $\mu_a(u)$ used for the definition of the refined knot vectors $\check U$ and $\tilde U$, respectively.}
	\label{linFun}
\end{figure}

To this end, and given a uniformly refined knot vector $\bar U$ with elements $\bar u_i$ (the analog derivation applies to $\bar V$ with elements $\bar v_i$), we shift the knots for $p_u, p_v >1$ (for $p_u, p_v =1$, we find that a uniform refinement strategy is suitable) linearly as
\begin{equation}
    \check{u}_i = \bar u_i + \Delta_{p_u} \mu(\bar u_i)
    \label{ucheck}
\end{equation}
with the linear function
\begin{equation}
    \mu(u) = \begin{cases}
       2u - 1 & \text{for} \quad u \in~]0,1[ \\[2mm]
       0 & \text{otherwise} 
    \end{cases}
\end{equation}
shown in Fig.~\ref{linFun} 
and the yet to be determined slope $\Delta_{p_u}$, which is fixed for the polynomial degree $p_u$.
The function $\mu(u)$ preserves the boundary knots at $u=0$ and $u=1$ (keeping the open knot vector structure intact) while shifting all internal knots by an amount proportional to their signed distance from the center of the parametric interval.
Corresponding to these knots $\check u_i$, the Greville sites of the refined basis can be expressed as
\begin{equation}
    \check\gamma_i^u = \bar\gamma_i^u + \dfrac{\Delta_{p_u}}{p_u-1}\sum_{j=1}^{p_u-1}\mu(\bar u_{i+j}) \,.
\end{equation}
The slope $\Delta_{p_u}$ is then determined by equating the $p_u$th Greville site of the primal and the dual basis, that is,
\begin{equation}
    \gamma_{p_u}^u = (\check\gamma_{2{p_u}-1}^u + \check\gamma_{2{p_u}}^u)/2
\end{equation}
which yields
\begin{equation}
    \Delta_{p_u} = \dfrac{2\gamma_{p_u}^u - \bar\gamma_{2{p_u}-1}^u  - \bar\gamma^u _{2{p_u}}}{{({p_u}-1)}^{-1} \sum_{j=1}^{{p_u}-1}\left[ \mu(\bar u_{2{p_u} - 1 + j}) + \mu(\bar u_{2{p_u} + j}) \right]} \,.
\end{equation}
This results in an exact match of the Greville sites of the primal and dual basis, except for the first and last $p_u-1$ Greville sites near the boundaries of the parametric domain.
The different behavior near the boundaries of the parametric domain originates from the decreasing distance between consecutive Greville sites only in this region, which is emphasized differently for the primal and refined basis due to the smaller knot spacing of the refined one. 
Correspondingly, a decreasing shift of the knots $\bar u_i$ is required as the parametric domain boundaries are approached, that is, a negative slope $\Delta_{p_u}$ in~\eqref{ucheck}.
In order to do so while maintaining a knot vector with well-separated knots (i.e., $\check u_{i+1} > \check u_i$ with a distance that is not arbitrarily close), we propose to ultimately modify the knot vector for the refined basis as
\begin{equation}
    \tilde{u}_i = \bar u_i + \Delta_{p_u} \mu_a(\bar u_i)
\end{equation}
with the clipped continuous function
\begin{equation}
    \mu_a(u) = \begin{cases}
        \dfrac{2a - 1}{a} u & \text{for}\quad  u \leq a \\[3mm]
       \, 2u - 1 & \text{for} \quad a < u \leq 1- a \\[2mm]
        \dfrac{2a - 1}{a} (u-1) & \text{for} \quad 1- a < u \leq 1
    \end{cases}
\end{equation}
also depicted in Fig.~\ref{linFun}, where we choose the parameter
\begin{equation}
    a = \dfrac{p_u}{2N_u - 1 - p_u} \,.
\end{equation}
The clipped function $\mu_a$ coincides with $\mu$ in the interior of the parametric domain but linearly transitions back to zero near the boundaries.
While this choice might not be optimal for the clipped region, we find that it is sufficient to obtain well-conditioned Gram matrices up to large values of $N_u$ and $N_v$, as demonstrated by the following numerical studies.

\begin{figure}
	\centering
	\includegraphics[scale=1]{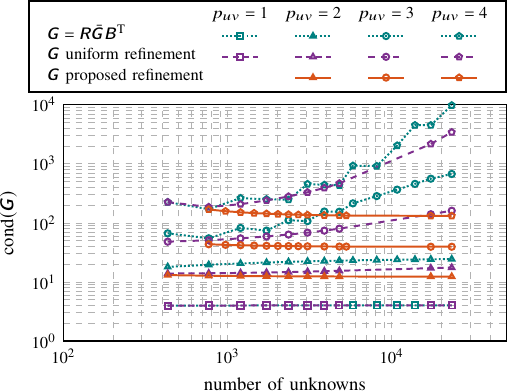}
	\caption{Conditioning of the Gram matrix for a closed surface for different polynomial degrees $p_u=p_v=p_{uv}$ and $N_{uv} \in \{ p_{uv}+1, \dots, 75 \}$.}
	\label{gram}
\end{figure}

\begin{figure}
	\centering
	\includegraphics[scale=1]{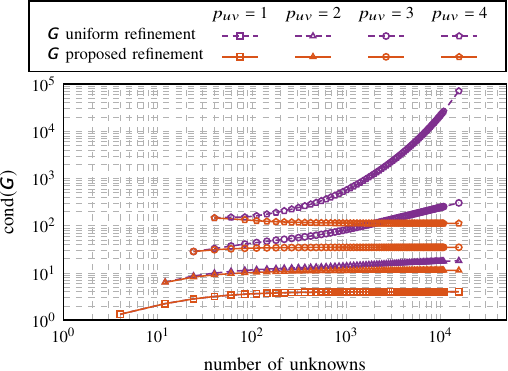}
	\caption{Conditioning of the Gram matrix for an open surface for different polynomial degrees $p_u=p_v=p_{uv}$ and $N_{uv} \in \{ p_{uv}+1, \dots, 75, 100 \}$.}
	\label{gramPlate}
\end{figure}

In Fig.~\ref{gram}, we have computed the condition numbers of the Gram matrix $\vec G$ for different polynomial degrees $p_u=p_v=p_{uv}$ and varied $N_u = N_v = N_{uv}$ for a closed object.
The same is done for an open object in Fig.~\ref{gramPlate}.
Note that the actual shape of the object is irrelevant, as the integral in~\eqref{gramDef} can be equivalently performed solely in the parametric domain, that is,
\begin{equation}
    {[\vec G]}_{mn} = \iint_{{[0,1]}^2} \big[\n \times \fhatddm~\big] \cdot \fhattilden  ~\,\dd u\dd v
    \label{Gparam}
\end{equation}
does neither depend on the Jacobi matrix $\veg J$ nor on its determinant $D$.
This can be shown independently of the specific choice of $\fhatddm$ and $\fhattilden$ from the properties of the covariant and contravariant curvilinear coordinates.
In both scenarios, the condition number of $\vec G$ remains bounded as $N_{uv}$ increases, independent of the polynomial degree.
This is in stark contrast to the uniformly refined case, where the condition number grows unboundedly with $N_{uv}$, and confirms the effectiveness of the proposed non-uniform refinement strategy.
For the latter, the obtained condition numbers are sufficiently small to allow for an efficient iterative inversion of the sparse $\vec G$ by means of a standard preconditioner such as an incomplete LU decomposition.
In Fig.~\ref{gram}, we have also included a refinement strategy where the primal and refined Bezier meshes are aligned, corresponding to the label $\vec R \bar{\vec G} \vec B^\T$.
Such a strategy might be desirable from an implementation perspective, as will be explained in Section~\ref{implement}, but it is, in general, not suitable to obtain a well-conditioned $\vec G$.

    \subsection{Dual Basis -- Properties} \label{secProperties}

Cross-checking the required properties i), ii), and iii) for the proposed primal and dual basis functions, we find that i) is clearly satisfied:
since the refined basis functions \smash{$\ghatijd$} are divergence-conforming by construction, any linear combination of them (including the dual basis functions) is divergence-conforming as well.
Property ii) of a well-conditioned Gram matrix $\vec G$ is achieved by the proposed non-uniform refinement strategy for the knot vectors $\bar U$ and $\bar V$, as detailed in the previous Section~\ref{secGramCond}.
Property iii), the vanishing of $\TPd \vec G^{-1} \TP$, is satisfied by associating the primal and dual basis functions dually to the edges of the Greville mesh $\mathscr{M}_\mr{G}$:
The roles of the loop and star basis mapping matrices $\matL$ and $\matS$ as defined in~\cite{hofmannLowFrequencyStabilizationBSpline2024} are interchanged for the dual basis.
As a consequence, the relations
\begin{equation}
	\TP\matL = \matO, \quad \matL^\T \TP = \matO, 
\end{equation}
\begin{equation}
	 \TP \vec H = \matO, \quad \vec H^\T \TP = \matO
\end{equation}
holding for the primal basis carry over to
\begin{equation}
	\TPd\matS = \matO, \quad \matS^\T \TPd = \matO, 
\end{equation}
\begin{equation}
	\TPd \widetilde{\vec H} = \matO, \quad \widetilde{\vec H}^\T \TPd = \matO
\end{equation}
for the dual basis, where $\vec H$ and $\widetilde{\vec H}$ denote the (potentially empty) primal and dual harmonic subspace mapping matrices for multiply-connected geometries.

With this, all three properties i), ii), and iii) have been established for the proposed primal and dual bases, and a suitable dual basis for the construction of a multiplicative Calderón preconditioner of arbitrary polynomial degree is at hand.
Moreover, in the lowest-order case, $p_u=p_v=1$, the primal basis reduces precisely to \ac{RWG} functions (also called rooftop or Raviart-Thomas functions) and the proposed dual basis to \ac{BC} functions on (curvilinear) quadrilateral subdomains of the patches.
This shows that a generalization of the \ac{BC} functions to arbitrary polynomial degrees has been obtained, retaining all of the advantageous properties known from the lowest-order case.

    \subsection{Dual Basis -- Implementation Aspects} \label{implement}

As pointed out in~\cite{andriulliMultiplicativeCalderonPreconditioner2008}, it might be desirable to express the matrices involved in the preconditioned system~\eqref{GTGTcalderon} as products of sparse mapping matrices and matrices assembled with respect to the refined basis $\ghatijd$.
Indeed, the discretization with dual basis functions of $\vecop T_{\!\tilde k}$ can be expressed as
\begin{equation}
    \widetilde{\vec T_{\!\tilde k}} = \vec B \bar{\vec T}_{\!\tilde k} \vec B^\T \,,
\end{equation}
where ${[\bar{\vec T}_{\!\tilde k}]}_{mn} = \iint_\Gamma \fbarm \cdot \big(\vecop T_{\!\tilde k}\,\fbarn  \big) \,\dd S(\veg r)$ is the discretization with the refined basis and $\vec B$ is a sparse matrix containing the coefficients $\alpha_{st}^{d}$.
However, constructing also the discretization with the primal basis functions from the refined one as
\begin{equation}
    \vec T_{\! k} = \vec R \bar{\vec T}_{\!k} \vec R^\T
    \label{viaR}
\end{equation}
and similarly
\begin{equation}
    \vec G = \vec R \bar{\vec G} \vec B^\T \,,
\end{equation}
is in general not possible.
The underlying reason is that the primal and refined Bezier meshes (defined by the knot vectors $U$ and $\bar U$ as well as $V$ and $\bar V$) are not necessarily aligned\footnote{The integration over the non-aligned Bezier meshes is only needed for the sparse Gram matrix $\vec G$ in~\eqref{Gparam}, where it can be performed at a negligible computational cost by forming an intersection mesh.}.
While they are aligned in the case $p_u=p_v=1$, they are in general not for $p_u>1$ and $p_v>1$.
In the case that they are aligned (or are deliberately chosen to be as done in Fig.~\ref{gram}), the sparse matrix $\vec R$ exists and can be constructed accurately via so-called knot insertion~\cite{pieglNURBSBook1995}.
However, since in many practical cases we have $k \neq \tilde k$ for the wavenumbers of the primal and dual discretizations anyway~\cite{adrianElectromagneticIntegralEquations2021,merliniMagneticCombinedField2020, DarbasDiss, boubendirWellconditionedBoundaryIntegral2014}, the non-existence of~\eqref{viaR} is only a minor drawback.
At the same time, the overall approach remains compatible with classical fast methods~\cite{zhaoAdaptiveCrossApproximation2005,chewFastEfficientAlgorithms2001,jukicStabilizedMultilevelBSplineBased2026}.

	\section{Numerical Results}

\begin{figure}[tp]
	\centering
	\subfloat[][NURBS]{\includegraphics[scale=1.0]{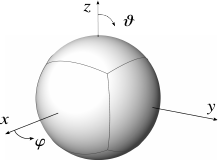}\label{geo}} 	\hfill
	\subfloat[][RWGs]{\includegraphics[scale=0.13]{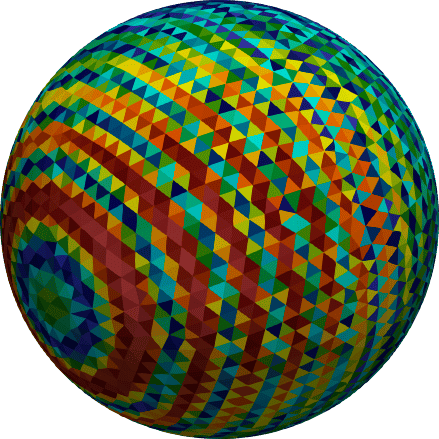}\label{rwg}} 	\hfill
	\subfloat[][$p_{uv}=1$]{\includegraphics[scale=0.13]{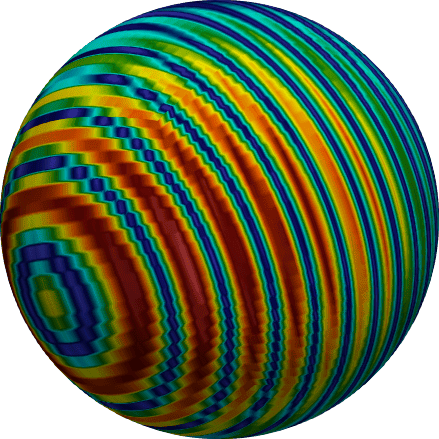}\label{p1}} \\
	\hfill \subfloat[][$p_{uv}=2$]{\includegraphics[scale=0.067]{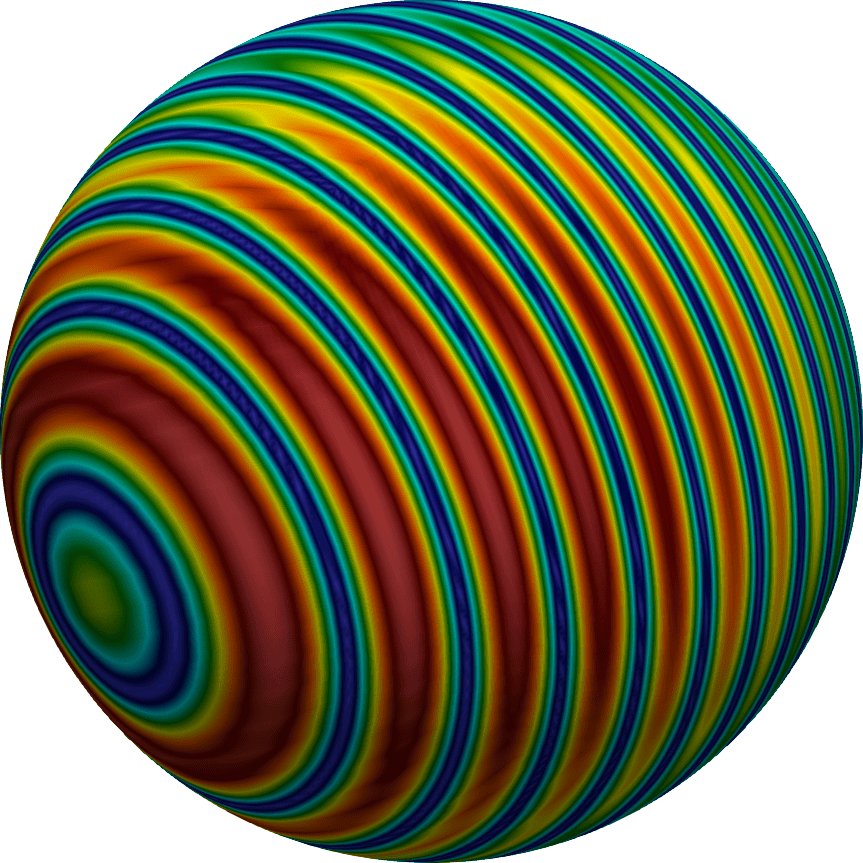}\label{p2}} \hspace{4mm}
	\subfloat[][$p_{uv}=3$]{\includegraphics[scale=0.067]{xBspline3.png}\label{p3}} 
	\caption{Scattering of a plane wave from a sphere of $r_\mr{s}=\SI{1}{\meter}$ at $f=\SI{1}{\giga\hertz}$: a) \ac{NURBS} description; b) real parts of the induced surface current density for a description by \num{8100} RWGs on a triangulation; c) by \num{8100} B-spline basis functions of $p_{uv}=1$, d) $p_{uv}=2$, and e) $p_{uv}=3$.}
	\label{sphereCurrent}
\end{figure} 
\begin{figure}[tp]
	\centering
	\captionsetup[subfloat]{labelformat=empty}
	\hfill\subfloat[][]{\includegraphics[scale=1]{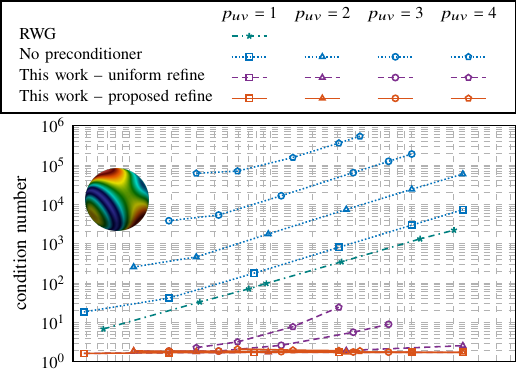}\label{sphereCond}} 	\\[-7mm]
	\hfill\subfloat[][]{\includegraphics[scale=1]{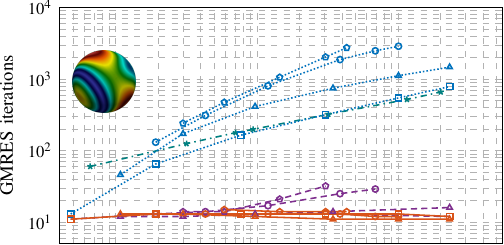}\label{sphereIter}} \\[-5mm]
	\hfill\subfloat[][]{\includegraphics[scale=1]{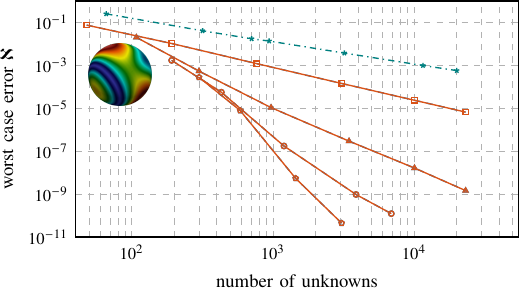}\label{sphereError}} 
	\caption{Scattering of a plane wave from a sphere of $r_\mr{s}=\SI{1}{\meter}$ at $f=\SI{100}{\mega\hertz}$ with and without proposed preconditioner.}
	\label{sphereVsN}
\end{figure}

To verify that the proposed dual basis yields an efficient Calderón preconditioner for arbitrary polynomial degrees, we consider the scattering of a plane wave from several \ac{PEC} objects.
In all of these, a \ac{GMRES} solver is employed without restarts and a relative residual of $\epsilon=\num{1e-12}$ as a stopping criterion (the residual is chosen this low due to a) the low error levels achieved in several of the examples and the need for fair comparisons and b) to study the behavior of the residual in more detail).
The implementation itself is based on~\cite{BEAST,nurbs2023} and the singularities in the numerical evaluation of the integrals in~\eqref{singular} and~\eqref{hsingular} are addressed by employing the regularizing coordinate transforms of~\cite[pp.~289ff]{sauterQuadratureHpGalerkin1997, Sauter2011}.
If not noted otherwise, we set $\tilde k = k$ for $\widetilde{\vec T}_{\!\tilde k}$ in the preconditioner~\eqref{GTGTcalderon}.

    \subsection{Scattering from a Sphere}

As a first scattering object, we consider a sphere of radius $r_\mr{s}=\SI{1}{\meter}$ described exactly (i.e., there is no approximation of the true sphere geometry) by 6 \ac{NURBS} patches as depicted in Fig.~\ref{sphereCurrent}~a).
Also shown in Fig.~\ref{sphereCurrent} are the real parts of the induced surface current densities for \num{8100} \ac{RWG} basis functions on a triangulation of the sphere and for \num{8100} basis functions as employed in this work for polynomial degrees of $p_{uv}=1,2,3$.
Clearly, the exact geometry description and the higher-order nature of the basis improve the solution accuracy.
To quantify the accuracy, we compute the error of the scattered fields at a distance of $r=\SI{5}{\meter}$ with respect to a Mie series expansion~\cite{jinTheoryComputationElectromagnetic2015,hofmannSphericalScatteringJuliaPackage2023}).
To do so, the fields are computed on a spherical grid with $\vartheta$ and $\varphi$ in steps of \SI{5}{\degree}, and the relative worst-case error
\begin{equation}
	\aleph = \max_{\vartheta, \varphi} \{   \frac{\left|e(\vartheta, \varphi) - \hat{\rule{0ex}{0.2ex}\mkern-3.0mu e}(\vartheta, \varphi)\right|}{\max\limits_{\vartheta, \varphi} \left|e(\vartheta, \varphi)\right|}    \} 
	\label{errorDef}
\end{equation}
is computed for the electric field. 
In all cases, we have also computed the errors for the magnetic field and the \ac{FF} and verified that they match the electric field error.
As can be seen from Fig.~\ref{sphereVsN}, the error converges as $h^{2p+1}$ with $h$ being the maximum distance between neighboring knots in $U$ and $V$, following theoretical guarantees~\cite{dolzIsogeometricBoundaryElements2019}.
Importantly, this is independent of the application of the proposed preconditioner, confirming the implementation's fundamental correctness.
Notably, compared to a triangulation approach employing \acp{RWG}, the lowest order B-spline-based basis functions (also constituting \acp{RWG} on quadrilateral subdomains) achieve for the same number of unknowns about two digits more in accuracy, presumably due to the absence of a geometrical approximation error.
From the condition numbers and the number of \ac{GMRES} iterations shown in Fig.~\ref{sphereVsN} for different number of unknowns $N$, we can observe: 
a) As expected, without preconditioner, the condition number and the iteration count increase with $N$ and the polynomial degree~\cite{valdesHighorderDivQuasi2011}.
b) The proposed preconditioner based on a uniform knot refinement strategy for the proposed dual basis significantly improves the condition number and iterations, but still leads to a growth in $N$.
c) With the proposed refinement strategy, the condition number and iterations become virtually independent of $N$ and the polynomial degree $p_{uv}$.
Specifically, a condition number of 1.9 and 11 iterations (for a residual of $\epsilon=\num{1e-12}$) are achieved.
This also highlights the importance and effectiveness of the refinement strategy underlying the dual basis.

\begin{figure}[tp]
	\centering
	\includegraphics[scale=1]{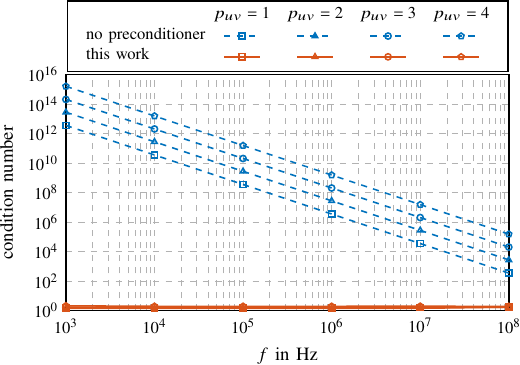}
	\caption{Condition number for a sphere of $r_\mr{s}=\SI{1}{\meter}$ discretized with \num{1452} unknowns for different frequencies.}
	\label{sphereVsFN12}
\end{figure}

Fixing the number of unknowns on the sphere to \num{1452} (which is possible due to the B-spline nature of the basis independent of $p_{uv}$) and decreasing the frequency, the condition number with and without the proposed preconditioning is depicted in Fig.~\ref{sphereVsFN12}.
As can be seen, regardless of the polynomial degree, the low-frequency breakdown is overcome, confirming the correctness of the proposed higher-order discretization of the Calderón preconditioner. 
As mentioned before, however, to fully overcome the breakdown, a combination with quasi-Helmholtz projectors is needed~\cite{andriulliWellConditionedElectricField2013,hofmannLowFrequencyStabilizationBSpline2024,hofmannExcitationAwareSelfAdaptiveFrequency2023,hofmannLowFrequencyStableDiscretization2022,hofmannLowFrequencyStableExcitation2023}.

    \subsection{Scattering from a Cube}

\begin{figure}[tp]
	\centering
	\captionsetup[subfloat]{labelformat=empty}
	\hfill\subfloat[][]{\includegraphics[scale=1]{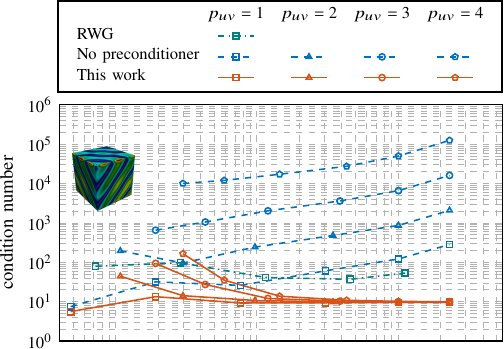}} 	\\[-7mm]
	\hfill\subfloat[][]{\includegraphics[scale=1]{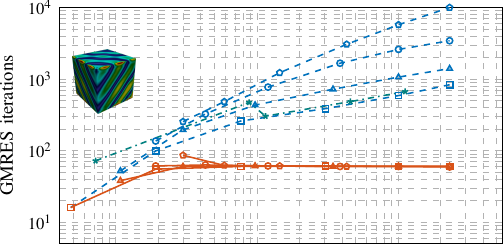}} \\[-5mm]
	\hfill\subfloat[][]{\includegraphics[scale=1]{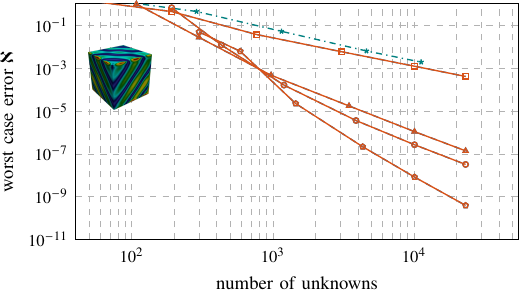}} 
	\caption{Scattering of a plane wave from a cube of edge length $a=\SI{1}{\meter}$ at $f=\SI{300}{\mega\hertz}$ with and without proposed preconditioner.}
	\label{cubeVsN}
\end{figure}

As a second example, we consider the scattering of a plane wave from a cube of edge length $a=\SI{1}{\meter}$ at a frequency of $f=\SI{300}{\mega\hertz}$ with the results depicted in Fig.~\ref{cubeVsN}.
The error in the scattered field is computed using the null-field condition inside the cube: it is measured by how well the scattered field cancels the incident field.
As can be seen, increasing the polynomial degree still improves the achieved error, but due to the edge and corner singularities, $h^{2p+1}$ is no longer satisfied.
Notably, even though an also shown \ac{RWG} discretization does not suffer from a geometrical approximation error, the employed basis seems to better use the degrees of freedom in the lowest order case (as also observed in~\cite{dolzNumericalComparisonIsogeometric2020}).

Regarding the preconditioner's effectiveness, the condition number is again significantly improved.
For a small number of unknowns, it is slightly increased but approaches a constant level of 10 already at around \num{1000} unknowns.
The number of iterations is constant at 22, all independent of the polynomial degree.

    \subsection{Scattering from an Open Geometry}

\begin{figure}[tp]
	\centering
	\captionsetup[subfloat]{labelformat=empty}
	\hfill\subfloat[][]{\includegraphics[scale=1]{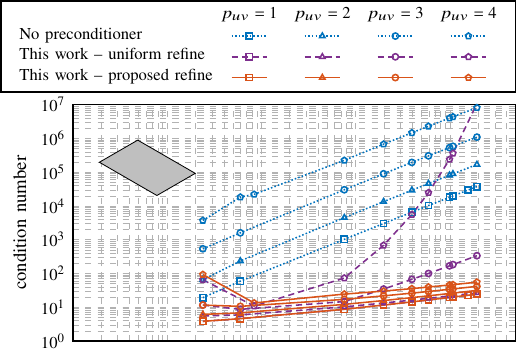}} 	\\[-7mm]
	\hfill\subfloat[][]{\includegraphics[scale=1]{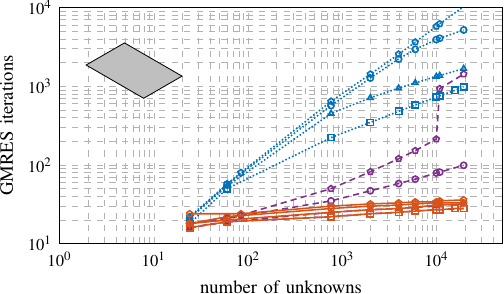}} 
	\caption{Scattering of a plane wave from a square plate of edge length $l=\SI{1}{\meter}$ at a frequency of $f=\SI{100}{\mega\hertz}$.}
	\label{plateVsN}
\end{figure}

To show that the proposed dual basis is also valid for open geometries, we consider the scattering from a square plate of edge length $l=\SI{1}{\meter}$ at a frequency of $f=\SI{100}{\mega\hertz}$ as depicted in Fig.~\ref{plateVsN}.
As with the scattering from the sphere, we see that the proposed non-uniform refinement strategy underlying the dual basis is crucial for obtaining an efficient preconditioner.
The slowly increasing condition number and iteration count are expected and also occur in the \ac{RWG} case as the Calderon identities do not strictly hold on screens~\cite{MathematicalModelsMethods}.
However, it still constitutes an efficient preconditioning strategy as corroborated by the low iteration counts below 36 for \num{19404} basis functions and $p_{uv}=4$.

    \subsection{Scattering from Realistic Geometries}

\begin{figure*}[tp]
	\centering
	\subfloat[][NURBS model]{\includegraphics[scale=0.18]{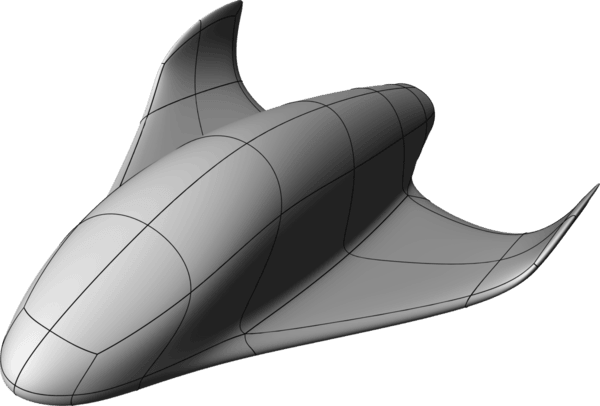}\label{geoShip}} 	\hfill
	\subfloat[][RWGs]{\includegraphics[scale=0.18]{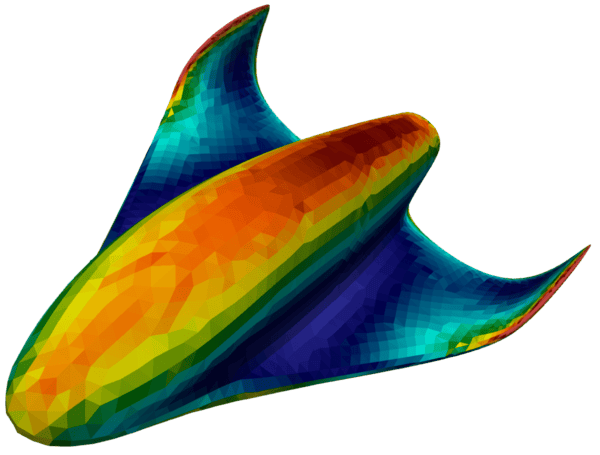}\label{rwgShip}} 	\hfill
	\subfloat[][$p_{uv}=1$]{\includegraphics[scale=0.18]{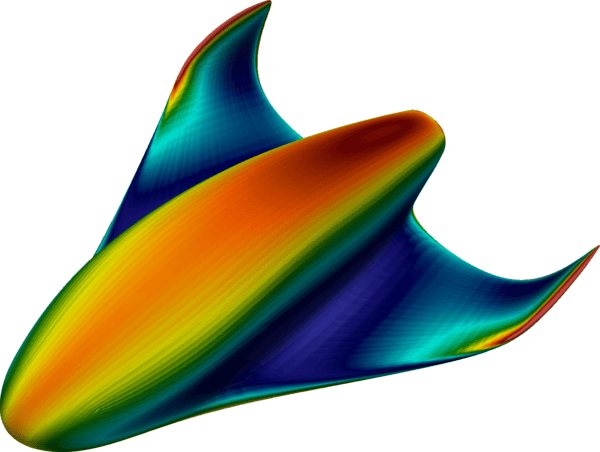}\label{p1Ship}} \hfill
	\subfloat[][$p_{uv}=2$]{\includegraphics[scale=0.18]{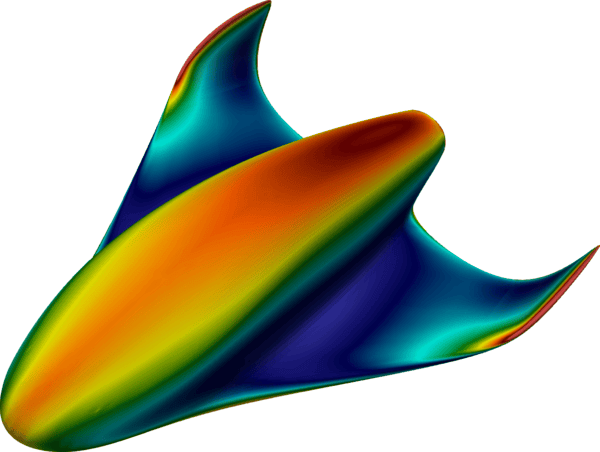}\label{p2Ship}} \\
	\captionsetup[subfloat]{labelformat=empty}
	\subfloat[][]{\includegraphics[scale=1]{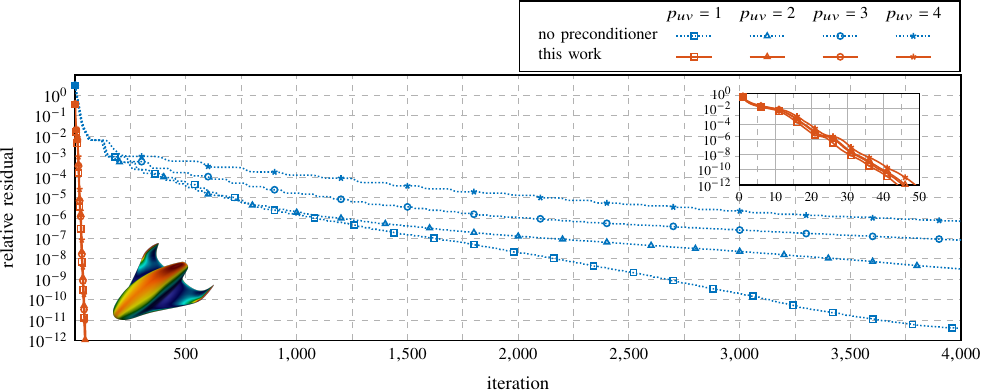}\label{p2ShipPlot}}
	\vspace{-3mm}
	\caption{Model of a spaceplane with \num{81} \ac{NURBS} patches, induced surface current densities and GMRES residual for the illumination with a plane wave, where the shuttle has a maximum extent of $0.08\lambda$ and a discretization with \num{23328} basis functions.}
	\label{spaceShipCurrent}
\end{figure*}

As a more realistic setup, the model of a spaceplane shown in Fig.~\ref{spaceShipCurrent}~a) is illuminated by a plane wave.
\begin{figure*}[tp]
	\centering
	\subfloat[][$p_{uv}=1$]{\includegraphics[scale=0.18]{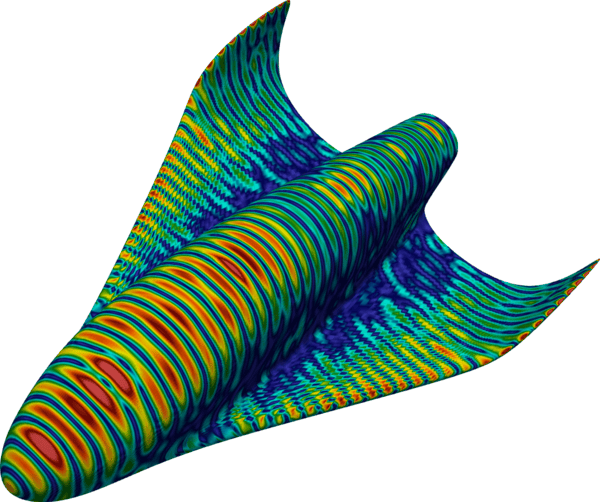}\label{spaceShipP1}} 	\hfill
	\subfloat[][$p_{uv}=2$]{\includegraphics[scale=0.18]{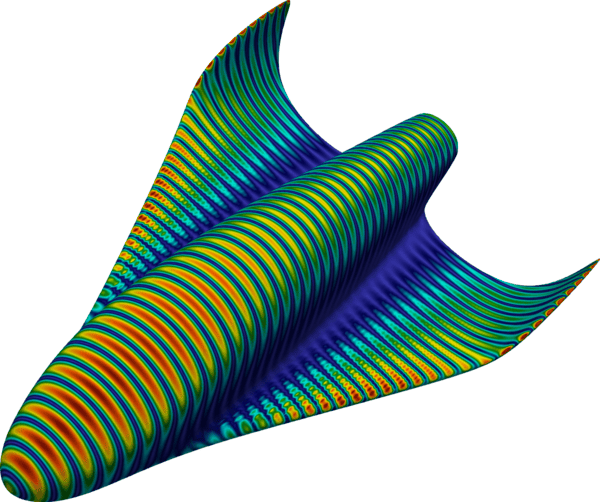}\label{spaceShipP2}} 	\hfill
	\subfloat[][$p_{uv}=3$]{\includegraphics[scale=0.18]{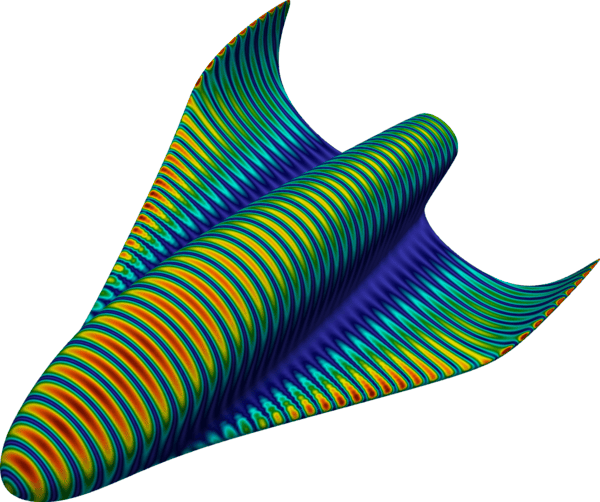}\label{spaceShipP3}}  	\hfill
	\subfloat[][mean curvature $H$]{\includegraphics[scale=1]{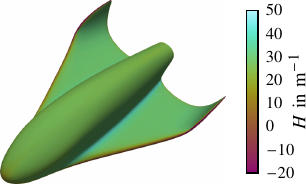}\label{spaceShipP4}} \\
	\captionsetup[subfloat]{labelformat=empty}
	\subfloat[][]{\includegraphics[scale=1]{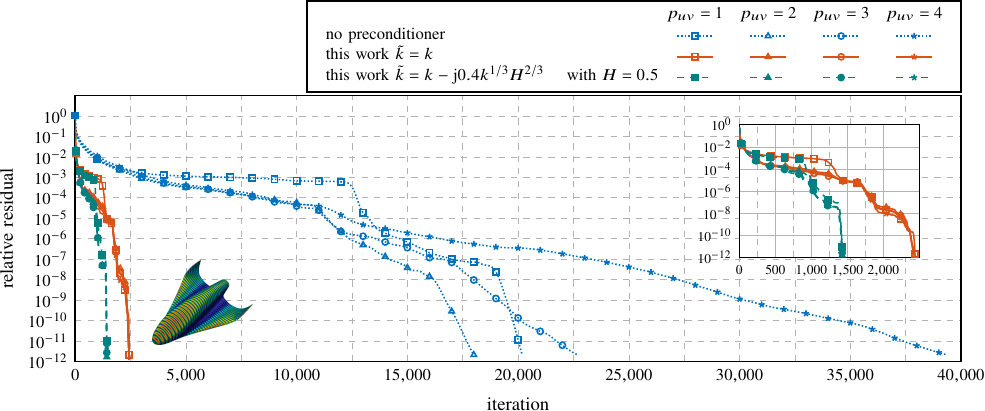}\label{spaceShipPlot}}
	\vspace{-3mm}
	\caption{Scattering from the spaceplane for the illumination with a plane wave, where the shuttle has a maximum extent of $21\lambda$ and is discretized with \num{109512} basis functions: induced surface current densities, mean curvature $H$ truncated to a maximum of \SI{50}{\per \meter}, and GMRES residual.}
	\label{spaceShip3GHz}
\end{figure*}
\begin{figure*}[tp]
	\centering
	\subfloat[][NURBS model]{\includegraphics[scale=0.32]{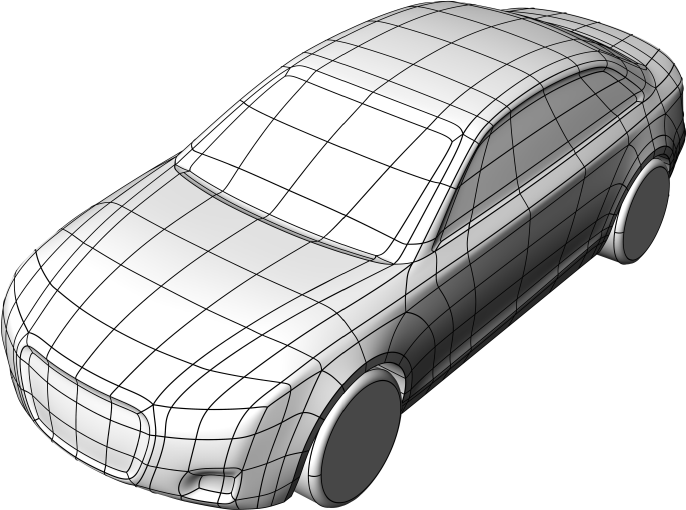}\label{geoCar}} 	\hfill
	\subfloat[][$p_{uv}=1$]{\includegraphics[scale=0.21]{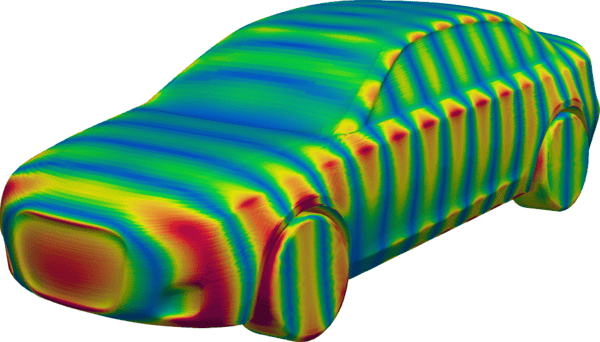}\label{currCar}} 	\hfill
	\subfloat[][$p_{uv}=2$]{\includegraphics[scale=0.21]{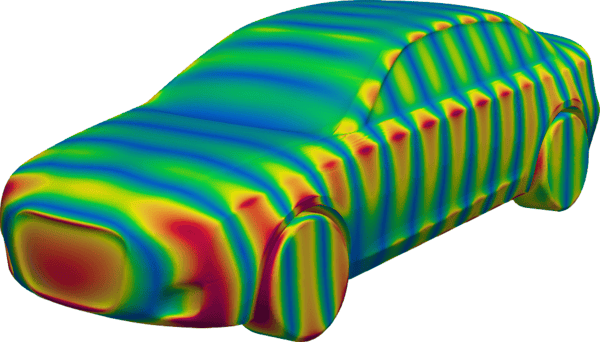}\label{currCar2}} \hfill
	\subfloat[][$p_{uv}=3$]{\includegraphics[scale=0.21]{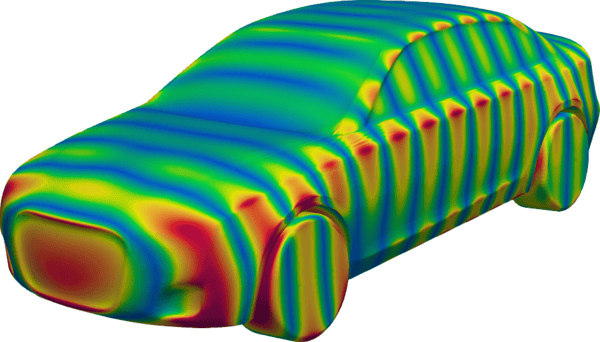}\label{currCar3}}  \\
	\captionsetup[subfloat]{labelformat=empty}
	\subfloat[][]{\includegraphics[scale=1]{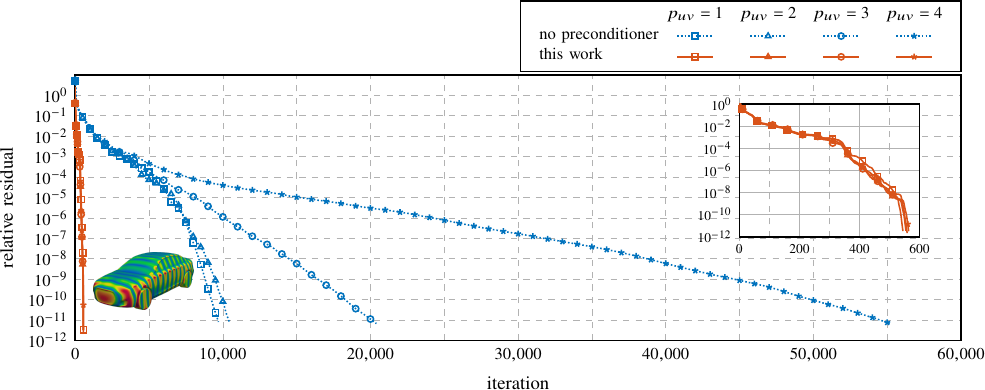}\label{currCarPlot}}	
	\vspace{-3mm}
	\caption{Model of a car with \num{1032} \ac{NURBS} patches, induced surface current densities and GMRES residual for the illumination with a plane wave, where the car has a maximum extent of $8.3\lambda$ and a discretization with \num{132072} basis functions.}
	\label{car}
\end{figure*}
The model consists of \num{81} \ac{NURBS} patches and is discretized with \num{23328} basis functions, identical for the polynomial degrees $p_{uv} = 1,2,3,4$.
In the first setup, we choose the frequency such that the shuttle has a maximum extent of $0.08\lambda$.
The fast convergence of \ac{GMRES} within 45 iterations, independent of the polynomial degree, compared to several thousand iterations (\num{3996} for $p_{uv} = 1$, \num{6952} for $p_{uv} = 2$, \num{12043} for $p_{uv} = 3$, \num{12488} for $p_{uv} = 4$) without preconditioning, demonstrates its effectiveness in this setup as well.

\begin{figure}[tp]
	\centering
	\includegraphics[scale=1]{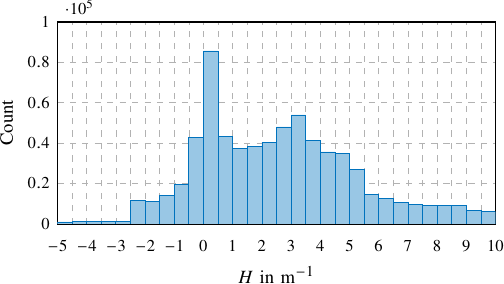}
	\caption{Extract of the histogram for mean curvature $H$ of the spaceplane showing the most often occurring values. The actual minimum occurring value is \SI{-20.06}{\per \meter} and the actual maximum occurring value \SI{262.2}{\per \meter}.}
	\label{spaceShuttleHist}
\end{figure}

\begin{figure}[tp]
	\centering
	\includegraphics[scale=1]{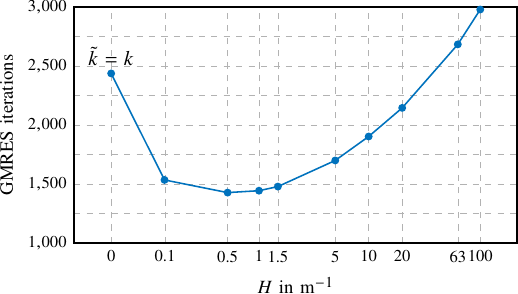}
	\caption{Number of iterations for the $21\lambda$ spaceplane for different values of $H$ in $\tilde k = k - \jm 0.4 k^{1/3}H^{2/3}$ for the Calderón preconditioner for $p_{uv}=2$.}
	\label{spaceShuttleIterVsH}
\end{figure}

Similarly, when increasing the frequency such that the spaceplane spans $21\lambda$, the induced current densities and \ac{GMRES} residual are depicted in Fig.~\ref{spaceShip3GHz}.
Accordingly, the number of basis functions is increased to \num{109512}.
Also, in this case, the preconditioner reduces the iteration count significantly: from up to \num{40000} to \num{2436}, independent of $p_{uv}$.
Since this example is in the high-frequency regime, choosing the wavenumber of the preconditioner as $\tilde k = k - \jm 0.4 k^{1/3}H^{2/3}$ can even reduce the number of iterations further~\cite{adrianElectromagneticIntegralEquations2021, DarbasDiss, boubendirWellconditionedBoundaryIntegral2014}.
This is confirmed by the achieved iteration count of \num{1427}, also shown in Fig.~\ref{spaceShip3GHz}.
However, while $H$ is often heuristically chosen to be the maximum mean curvature of the scatterer's surface, we find that a lower value results in fewer iterations.
Specifically, we have computed the mean curvature for the spaceplane (evaluated on a 101 $\times$ 101 point grid on each patch via~\cite{nurbs2023}) and depicted it in Fig.~\ref{spaceShip3GHz}~d) with a truncation to a maximum of \SI{50}{\per \meter}.
The actual maximum is about \SI{262}{\per \meter}.
As can be seen, the mean curvature is relatively low on most of the surface and only has high values on the tips of the wings.
More precisely, a histogram for the mean curvature in Fig.~\ref{spaceShuttleHist} shows that on most of the surface, the mean curvature is between \SI{0}{\per \meter} and \SI{0.5}{\per \meter}. 
Accordingly, choosing $H$ for the computation of $\tilde k$ as \SI{0.5}{\per \meter} (i.e., the mode of samples of $H$) resulted in the lowest number of iterations as the study in Fig.~\ref{spaceShuttleIterVsH} shows.
The solution times for the setups without preconditioner, with preconditioner setting $\tilde k = k$, and setting $\tilde k = k - \jm 0.4 k^{1/3}H^{2/3}$ are summarized in Table~\ref{timingsShuttle}.
\begin{table}[tp]
    \centering
    \caption{Solution times for the scattering from the $21 \lambda$ spaceplane.}
    \label{timingsShuttle}
    \begin{tabular}{lcccc}
        \toprule
         & $p_{uv}=1$ & $p_{uv}=2$ & $p_{uv}=3$ & $p_{uv}=4$ \\
        \cmidrule(lr){2-2}\cmidrule(lr){3-3}\cmidrule(lr){4-4}\cmidrule(lr){5-5}
        No preconditioner & 5.7 h & 8.3 h & 10.6 h & 19.3 h \\
        This work $\tilde k = k$   & 54 min & 54 min & 54 min & 54 min \\
        This work $\tilde k = k(H)$   & 20 min & 20 min & 20 min & 20 min \\
        \bottomrule
    \end{tabular}
\end{table}
Note that we did not employ any acceleration techniques, but explicitly constructed the matrices and performed the full matrix-vector products utilizing \num{512} threads on a sysGen Server G493-ZB1-AAP1 equipped with two dual AMD EPYC TURIN 9755, each with 128 cores and a clock frequency of $\SIrange{2.70}{4.10}{GHz}$, and $\SI{6}{TB}$ DDR5 RAM.

As a last example, we consider the model of a car described by \num{1032} \ac{NURBS} patches and \num{132072} basis functions as shown in Fig.~\ref{car}, where 
the iteration count is reduced from up to \num{55000} to about \num{550}, again independent of the polynomial degree.
The solution times with and without the proposed preconditioner are summarized in Table~\ref{timingsCar} showing the effectiveness of the preconditioner.
\begin{table}[tp]
    \centering
    \caption{Solution times for the scattering from the car.}
    \label{timingsCar}
    \begin{tabular}{lcccc}
        \toprule
         & $p_{uv}=1$ & $p_{uv}=2$ & $p_{uv}=3$ & $p_{uv}=4$ \\
        \cmidrule(lr){2-2}\cmidrule(lr){3-3}\cmidrule(lr){4-4}\cmidrule(lr){5-5}
        No preconditioner & 2.3 h & 2.75 h & 11.3 h & 37.0 h \\
        This work         & 17.8 min & 17.8 min & 17.8 min & 17.8 min \\
        \bottomrule
    \end{tabular}
\end{table}

    \section{Conclusion}

We have proposed the first explicit higher-order dual basis suitable for multiplicatively Calderón preconditioning the \ac{EFIE}.
It can be considered a direct generalization of the \ac{BC} functions to arbitrary polynomial degrees and constitutes a fundamental building block for discretizations and preconditioning strategies relying on a dual basis.
The underlying B-spline-based basis functions offer more flexibility due to their local support and the independence of the polynomial degree from the number of splines on a patch.

The obtained preconditioner has been shown to significantly reduce the number of iterations required to solve the \ac{EFIE}.
Most notably, in the considered examples, the iteration count remained fixed, independent of the polynomial degree, which is key to fully leverage the potential of higher-order bases and to overcome long solution times caused by increased condition numbers for higher polynomial degrees.


	
	%

%

%

	\ifCLASSOPTIONcaptionsoff
	  \newpage
	\fi

	

	\printbibliography

	
	

\end{document}

%% file: praemble/hft-paper-praemble-main.tex
\input{hft-paper-praemble-packages.tex}
\input{hft-macros.tex}

\input{sba-macros.tex}

%% file: praemble/hft-paper-praemble-packages.tex

\usepackage[utf8]{luainputenc}

\usepackage[T1]{fontenc}

\usepackage{graphicx}
\usepackage[dvipsnames]{xcolor}

\usepackage[english]{babel}
\addto\captionsenglish{\renewcommand{\figurename}{Fig.}}
\addto\captionsenglish{\renewcommand{\tablename}{Table}}
\usepackage{csquotes}

\usepackage{newtxtext}
\usepackage{amsthm}
\usepackage[slantedGreek]{newtxmath}
\usepackage[OMLmathsfit]{isomath}
\DeclareMathAlphabet{\mathbfsf}{\encodingdefault}{\sfdefault}{bx}{n}
\usepackage{bm}
\usepackage{envmath}
\usepackage{mathtools}
\usepackage{commath}
\usepackage{siunitx}
\usepackage{nicefrac}

\usepackage[caption=false,font=footnotesize]{subfig}

\usepackage{booktabs}
\usepackage{footmisc}  

\usepackage{url}

\input{hft-packages-theorem.tex}

\input{hft-packages-review.tex}
\input{hft-packages-tikz-pgfplots.tex}
\usepackage[hidelinks]{hyperref} 
\input{hft-packages-cleveref.tex}

%% file: praemble/hft-packages-theorem.tex
\theoremstyle{definition}

\theoremstyle{plain}

\theoremstyle{remark}

%% file: praemble/hft-packages-review.tex
\usepackage{lineno}
\modulolinenumbers[5]
\usepackage{todonotes}
\usepackage{umoline}


%% file: praemble/hft-packages-tikz-pgfplots.tex
\usepackage{pgfplots}
\usepackage{pgfplotstable}
\pgfplotsset{compat=newest}
\pgfplotsset{plot coordinates/math parser=false}
\newlength\figureheight
\newlength\figurewidth
\pgfplotsset{every axis plot/.append style={line width=1.5pt},
    legend style={font=\footnotesize, 
        text height=1.0ex,
        draw=black,
        fill=white,
        legend cell align=left}}

%% file: praemble/hft-packages-cleveref.tex
\usepackage[english]{cleveref}

\Crefname{defn}{definition}{definitions}
\Crefname{defn}{Definition}{Definitions}

\Crefname{asm}{assumption}{assumptions}
\Crefname{asm}{Assumption}{Assumptions}

\crefname{lem}{lemma}{lemmas} 
\Crefname{lem}{Lemma}{Lemmas}

\crefname{prop}{proposition}{propositions} 
\Crefname{prop}{Proposition}{Propositions}

\crefname{thm}{theorem}{theorms} 
\Crefname{thm}{Theorem}{Theorms}

\crefname{cor}{corollary}{corollaries}
\Crefname{cor}{Corollary}{Corollaries}

%% file: praemble/hft-macros.tex
\input{hft-macros-environments.tex} 
\input{hft-macros-maths.tex} 
\input{hft-packages-acro.tex} 
\input{hft-macros-booktabs.tex}


%% file: praemble/hft-macros-environments.tex
\newcounter{subequation}
\newlength\mtabskip\mtabskip=-1.25cm

\def\mtabLong{long}
\makeatletter

\makeatother

%% file: praemble/hft-macros-maths.tex
\newcommand{\mr}{\mathrm}

\newcommand{\veg}[1]{\bm{#1}}     
\newcommand{\mat}[1]{\mathsfbfit{#1}} 
\renewcommand{\vec}[1]{\mathsfbfit{#1}} 
\newcommand{\vecop}[1]{\bm{\mathcal{#1}}} 


\newcommand{\matO}{\mathbfsf{0}}

\newcommand{\n}{\hat{\bm{n}}}


\newcommand{\dd}{\mathrm{d}}  


\newcommand{\jm}{\mathrm{j}}  

\newcommand{\e}{\mathrm{e}}





\newcommand{\T}{\mr{T}}


\newcommand\restr[2]{{
        \left.\kern-\nulldelimiterspace 
        #1 
        \vphantom{|} 
        \right|_{#2} 
}}

\newcommand\rst[3]{{
        \left.\kern-\nulldelimiterspace 
        #1 
        \vphantom{|} 
        \right|_{#2}^{#3} 
}}



%


\newcommand{\TAop}{\vecop T_{\!\!\mr{A},k}}
\newcommand{\TA}{\vec T_{\!\!\mr{A},k}}

\newcommand{\TPop}{\vecop T_{\!\!\upPhi,k}}
\newcommand{\TP}{\vec T_{\!\upPhi,k}}

\newcommand{\TPd}{\widetilde{\vec T}_{\!\upPhi,\tilde k}}

\newcommand{\fhatij}[0]{\,\,\hat{\mskip-6.0mu\veg f}_{\!ij}} 
\newcommand{\fhatn}[0]{\,\,\hat{\mskip-6.0mu\veg f}_{\!\!n}}

\newcommand{\ftil}[0]{\,\,\widetilde{\mskip-6.0mu\veg f}} 
\newcommand{\ftilij}[0]{\,\,\widetilde{\mskip-6.0mu\veg f}_{\!ij}}

\newcommand{\fhatddm}[0]{\,\,\hat{\mskip-6.0mu\veg f}_{\!\!m}} 
\newcommand{\fhatddn}[0]{\,\,\hat{\mskip-6.0mu\veg f}_{\!\!n}} 



\newcommand{\fhattilden}[0]{\fhatddn\hspace{-2.45mm}\overset{\sim}{\phantom{f}}}

\newcommand{\fhatqq}[0]{\,\,\hat{\mskip-6.0mu\veg f}_{\!\!ij}} 
\newcommand{\fhatqqu}[0]{\fhatqq{\vphantom{\veg f}}^{\!\!\!\!u}}
\newcommand{\fhattildeuij}[0]{\fhatqqu\hspace{-3.25mm}\overset{\sim}{\phantom{f}}} 

\newcommand{\fhatqqv}[0]{\fhatqq{\vphantom{\veg f}}^{\!\!\!\!v}}
\newcommand{\fhattildevij}[0]{\fhatqqv\hspace{-3.25mm}\overset{\sim}{\phantom{f}}} 

\newcommand{\fhatijd}[0]{\fhatij{\vphantom{\veg f}}^{\!\!\!\!d}}

\newcommand{\fhatqqd}[0]{\fhatqq{\vphantom{\veg f}}^{\!\!\!\!d}}

\newcommand{\ghatijd}[0]{\fhatqqd\hspace{-2.7mm}\smash{\raisebox{9.5pt}{\rule{0.8ex}{0.5pt}}}~\,}
\newcommand{\ghatiju}[0]{\fhatqqu\hspace{-2.7mm}\smash{\raisebox{10.5pt}{\rule{0.8ex}{0.5pt}}}~\,}
\newcommand{\ghatijv}[0]{\fhatqqv\hspace{-2.7mm}\smash{\raisebox{10.5pt}{\rule{0.8ex}{0.5pt}}}~\,}

\newcommand{\fhatst}[0]{\,\,\hat{\mskip-6.0mu\veg f}_{\!\!st}} 
\newcommand{\fhatstu}[0]{\fhatst{\vphantom{\veg f}}^{\!\!\!\!u}}
\newcommand{\fhatstv}[0]{\fhatst{\vphantom{\veg f}}^{\!\!\!\!v}}

\newcommand{\ghatstu}[0]{\fhatstu\hspace{-2.7mm}\smash{\raisebox{10.5pt}{\rule{0.8ex}{0.5pt}}}~\,}
\newcommand{\ghatstv}[0]{\fhatstv\hspace{-2.7mm}\smash{\raisebox{10.5pt}{\rule{0.8ex}{0.5pt}}}~\,}

\newcommand{\fbarm}[0]{\,\,\bar{\mskip-6.0mu\veg f}_{\!\!m}} 
\newcommand{\fbarn}[0]{\,\,\bar{\mskip-6.0mu\veg f}_{\!\!n}} 

%% file: praemble/hft-packages-acro.tex
\usepackage{acro}

\DeclareAcronym{DG}
{
    short = DG ,
    long = discontinuous Galerkin
}

\DeclareAcronym{ACA}
{
    short = ACA ,
    long = adaptive cross approximation
}

\DeclareAcronym{EFIE}
{
    short =  EFIE ,
    long = electric field integral equation
}

\DeclareAcronym{MFIE}
{
    short =  MFIE ,
    long = magnetic field integral equation
}

\DeclareAcronym{CFIE}
{
    short =  CFIE ,
    long = combined field integral equation
}

\DeclareAcronym{EFIO}
{
    short =  EFIO ,
    long = electric field integral operator
}

\DeclareAcronym{MFIO}
{
    short =  MFIO ,
    long = magnetic field integral operator
}

\DeclareAcronym{CFIO}
{
    short =  CFIO ,
    long = combined field integral operator
}

\DeclareAcronym{MUIE}
{
    short =  MUIE ,
    long = Müller integral equation
}

\DeclareAcronym{PMCHWT}
{
    short =  PMCHWT ,
    long = Poggio-Miller-Chang-Harrington-Wu-Tsai integral equation
}

\DeclareAcronym{SPD}
{
    short =  SPD ,
    long = {symmetric, positive definite}
}

\DeclareAcronym{SPSD}
{
    short =  SPD ,
    long = {symmetric, positive semi-definite}
}

\DeclareAcronym{PEC}
{
    short =  PEC ,
    long = perfectly electrically conducting
}

\DeclareAcronym{RWG}
{
    short = RWG ,
    long = Rao-Wilton-Glisson
} 

\DeclareAcronym{BC}
{
    short = BC ,
    long = Buffa-Christiansen
}

\DeclareAcronym{SVD}
{
    short = SVD ,
    long = singular value decomposition
}

\DeclareAcronym{CG}
{
    short = CG ,
    long = conjugate gradient
} 

\DeclareAcronym{PCG}
{
    short = PCG ,
    long = preconditioned conjugate gradient
} 

\DeclareAcronym{CGS}
{
    short = CGS ,
    long = conjugate gradient squared
}

\DeclareAcronym{CMP}
{
    short = CMP ,
    long = Calderón multiplicative preconditioner
} 

\DeclareAcronym{RFCMP}
{
    short = RF-CMP ,
    long = refinement-free Calderón multiplicative preconditioner
} 

\DeclareAcronym{HPD}
{
    short = HPD ,
    long = {Hermitian, positive definite}
} 

\DeclareAcronym{RHS}
{
    short = RHS ,
    long = right-hand side
}

\DeclareAcronym{LSE}
{
    short = LSE ,
    long = linear system of equations
}

\DeclareAcronym{AMG}
{
    short = AMG ,
    long = algebraic multigrid
}

\DeclareAcronym{PW}
{
    short = PW ,
    long = plane wave
}

\DeclareAcronym{GMRES}
{
    short = GMRES ,
    long = generalized minimum residual
}

\DeclareAcronym{IDR}
{
    short = IDR ,
    long = induced dimension reduction
}

\DeclareAcronym{BICGstab}
{
    short = BiCGstab ,
    long = stabilized bi-conjugate gradient
}

\DeclareAcronym{FF}
{
    short = FF ,
    long = far field
}

\DeclareAcronym{NF}
{
    short = NF ,
    long = near field
}

\DeclareAcronym{TE}
{
    short = TE,
    long = transverse electric
}

\DeclareAcronym{TM}
{
    short = TM,
    long = transverse magnetic
}

\DeclareAcronym{NURBS}
{
    short = NURBS,
    long = non-uniform rational B-splines
}

\DeclareAcronym{FEM}
{
    short = FEM,
    long = finite element method
}

\DeclareAcronym{PEEC}
{
    short = PEEC,
    long = partial element equivalent circuit
}

\DeclareAcronym{GWP}
{
    short = GWP,
    long = Graglia Wilton Peterson
}

%% file: praemble/hft-macros-booktabs.tex

\newcolumntype {n}{c}
\newcolumntype {N}{>{\small}c}
\newcolumntype {L}{>{\small}l}
\newcolumntype {F}{>{\footnotesize}c}
\newcolumntype {v}[1]{>{\raggedright \hspace {0pt}} p {#1}}
\newcolumntype {V}[1]{>{\small \raggedright \hspace {0pt}} p {#1}}
\newcolumntype{d}[1]{>{\DC@{.}{.}{#1}}c<{\DC@end}}

%
\newcolumntype{R}[1]{%
    >{\begin{turn}{90}\begin{minipage}{#1}\small\raggedright\hspace{0pt}}l%
            <{\end{minipage}\end{turn}}%
}

%% file: praemble/sba-macros.tex


\newcommand{\matL}{\mat{\Lambda}}   
\newcommand{\matS}{\mat{\Sigma}}